\documentclass[11pt]{article}
\usepackage[a4paper,margin=1in]{geometry}
\usepackage{amsmath,amssymb,amsthm}
\usepackage{booktabs}
\usepackage{array}
\usepackage{xcolor}
\usepackage{needspace}
\usepackage{cite}
\usepackage[colorlinks=true,linkcolor=blue!50!black,citecolor=blue!50!black]{hyperref}

\newtheorem{theorem}{Theorem}[section]
\newtheorem{lemma}[theorem]{Lemma}
\newtheorem{proposition}[theorem]{Proposition}
\newtheorem{corollary}[theorem]{Corollary}
\theoremstyle{remark}
\newtheorem{remark}[theorem]{Remark}

\newcommand{\GP}{\mathrm{GP}}

\title{\texorpdfstring{Perfect $3$-colorings of the generalized Petersen graphs $\GP(n,2)$}{Perfect 3-colorings of the generalized Petersen graphs GP(n,2)}}
\author{Hamed Karami\\[4pt]
\small Department of Mathematics and Statistics, Georgia State University, Atlanta, GA, USA\\
\small \texttt{hkarami1@student.gsu.edu}}
\date{}

\begin{document}
\maketitle

\begin{abstract}
We enumerate the parameter matrices of all perfect $3$-colorings of the generalized Petersen
graphs $\GP(n,2)$, $n\ge 5$. We show that there are only five matrices that can be the
parameter matrix of such a coloring. Two of them occur if and only if $5\mid n$, two occur if and only if $3\mid n$, and the remaining
one occurs if and only if $10\mid n$.
\end{abstract}



\section{Introduction}

Equitable partitions of graphs, also known as perfect colorings, arise in many parts of mathematics
and its applications; see \cite[Section~9.3]{GR}, \cite{God97}, and the survey \cite{KP25}. In coding
theory they are closely tied to completely regular codes, which go back to Delsarte \cite{Del73} and
include the perfect codes; see \cite{Neu92,BRZ19,KP25}. Delsarte's conjecture that Johnson graphs have no nontrivial perfect codes
\cite{Del73} has motivated the study of perfect $2$-colorings of Johnson graphs \cite{AMog11,GG13}.
Fon-Der-Flaass proved a bound on the correlation immunity of unbalanced Boolean functions and used it
to obtain a necessary condition for perfect colorings of hypercubes \cite{FDFb,FDFh}. In finite
geometry, intriguing and tight sets of point-line geometries fit into the same framework
\cite{DBS10}.

Equitable partitions are equally important outside combinatorics. The color refinement procedure,
also known as the one-dimensional Weisfeiler--Leman algorithm, computes the coarsest equitable
partition of a graph. It is a basic step of practical graph isomorphism software \cite{MP14}; two
graphs are fractionally isomorphic exactly when they have a common coarsest equitable partition
\cite{RSU94}; and the same partition limits what message-passing graph neural networks can
distinguish \cite{XHLJ19,Mor19}. In networks of identical dynamical units, equitable partitions
describe the patterns of synchrony that are forced by the network structure alone \cite{SGP03};
together with symmetries they organize cluster synchronization in networks of oscillators
\cite{Pec14,Sch16}, and they are used to characterize the controllability of leader--follower
consensus networks \cite{RJME09}. The related almost equitable partitions describe Laplacian
eigenvalues and eigenvectors \cite{CDR07}. Finally, for a $k$-regular graph an equitable partition
makes the simple random walk lumpable \cite{KS60}, so the walk can be studied through a much smaller
Markov chain on the parts.

For a fixed graph or family of graphs, a central problem is to determine all parameter matrices of
perfect colorings with a given small number of colors. Besides the hypercubes and Johnson graphs
mentioned above, perfect colorings of infinite circulant graphs with sufficiently many colors have
been described completely \cite{LA24}. Cubic graphs form a natural testing ground, since their
parameter matrices have row sums three and there are only a few candidates for each number of
colors. Avgustinovich and Lisitsyna \cite{AL11} described the perfect $2$-colorings of all transitive
cubic graphs with at most $18$ vertices, and the perfect colorings of prisms and M\"obius ladders are
known \cite{Lis,LA16}. Perfect $3$-colorings have been determined for several symmetric cubic graphs \cite{AQ}.

The generalized Petersen graphs $\GP(n,k)$ are a classical family of cubic graphs, introduced by
Coxeter \cite{Cox50} and named by Watkins \cite{Wat69}. The family contains the prisms
$\GP(n,1)$, the Petersen graph $\GP(5,2)$, the M\"obius--Kantor graph $\GP(8,3)$, the dodecahedron
$\GP(10,2)$, and the Desargues graph $\GP(10,3)$, and it has been studied from many points of view,
including automorphism groups \cite{FGW71}, Hamiltonicity \cite{Als83}, and perfect and total perfect
codes \cite{EJM09,WZ}. For two colors, perfect colorings of generalized Petersen graphs are well
understood: Alaeiyan and Karami \cite{AK} enumerated the parameter matrices of perfect $2$-colorings
of $\GP(n,2)$, Karami \cite{Kar} treated $\GP(n,3)$, and the case of arbitrary $k$ was treated in
\cite{Kar26} and, in the equivalent language of regular sets, by Wang and Zhang \cite{WZr}. The
perfect codes and total perfect codes mentioned above correspond to two of the possible parameter
matrices of perfect $2$-colorings, and for them the existence conditions were obtained in
\cite{EJM09} and \cite{WZ}, respectively. For three colors much less is known. Alaeiyan, Karami,
and Siasat \cite{AKS} determined
the parameter matrices for $\GP(5,2)$ and $\GP(6,2)$ and proved that $\GP(7,2)$ admits no perfect
$3$-coloring, Alaeiyan and Shokoohi \cite{AS} showed that $\GP(8,2)$ admits none, and the Petersen
graph and the dodecahedron are also covered by \cite{AQ}. Thus, for $\GP(n,2)$ only the
values $n=5,6,7,8,10$ had been settled, and to our knowledge the prisms are the only infinite family
of generalized Petersen graphs whose perfect $3$-colorings were completely known.

In this paper we classify all parameter matrices of perfect $3$-colorings of the generalized
Petersen graphs $\GP(n,2)$, $n\ge 5$, and determine for each of them the values of $n$ for which
it occurs.

\section{Preliminaries}\label{sec:prelim}

Let $G$ be a graph and let $m\ge 1$. A \emph{perfect $m$-coloring} of $G$ is a surjective
map $T\colon V(G)\to\{1,\dots,m\}$ for which there are integers $a_{ij}$ such that every vertex
of color $i$ has exactly $a_{ij}$ neighbors of color $j$, for all $i,j\in\{1,\dots,m\}$. The matrix
$A=(a_{ij})$ is its \emph{parameter matrix}.

Let $n\ge 3$ and $1\le k<n/2$. The \emph{generalized Petersen graph} $\GP(n,k)$ has vertices $a_i$
and $b_i$ for $i\in\{0,1,\dots,n-1\}$, and edges $a_ia_{i+1}$, $a_ib_i$, and $b_ib_{i+k}$, where
indices are always taken modulo $n$; it is a connected cubic graph.

For a perfect $3$-coloring of a connected cubic graph, let $v_i>0$ be the size of color class $i$.
Counting the edges between two classes in two ways gives
\begin{equation}\label{eq:balance}
v_ia_{ij}=v_ja_{ji}.
\end{equation}
Consequently, $a_{ij}=0$ if and only if $a_{ji}=0$, and multiplying the three instances of
\eqref{eq:balance} for the pairs $(1,2)$, $(2,3)$, $(3,1)$ and cancelling $v_1v_2v_3$ gives
\begin{equation}\label{eq:cycle}
a_{12}a_{23}a_{31}=a_{13}a_{32}a_{21}.
\end{equation}
See also~\cite[Theorem~2.1(2)]{DF21}. Every row sums to three. Consider the \emph{graph of color classes}, whose vertices are the three
color classes, two classes being adjacent when $a_{ij}>0$. Since the colored graph is connected, this
graph is connected; as it has three vertices, it is a path or a triangle.

\begin{lemma}\label{lem:candidates}
Every parameter matrix of a perfect $3$-coloring of a connected cubic graph is, up to renaming
colors, one of fifteen path matrices or three triangle matrices. The path matrices have the form
\begin{equation}\label{eq:path}
A(p,s,q,r)=\begin{pmatrix}3-p&p&0\\ q&3-q-r&r\\ 0&s&3-s\end{pmatrix},
\end{equation}
where either $q=r=1$ and $1\le p\le s\le 3$, or $(q,r)=(1,2)$ and $p,s\in\{1,2,3\}$. The triangle
matrices are
\[
A_3=\begin{pmatrix}1&1&1\\1&1&1\\1&1&1\end{pmatrix},\qquad
A_4=\begin{pmatrix}0&1&2\\1&1&1\\2&1&0\end{pmatrix},\qquad
C=\begin{pmatrix}0&1&2\\1&0&2\\1&1&1\end{pmatrix}.
\]
\end{lemma}

\begin{proof}
The same eighteen matrices are listed in \cite[Section~2]{AQ}; a longer list of twenty-two
candidates appears in \cite[Section~3]{AKS}, which also includes four matrices violating
\eqref{eq:cycle}. We include the short derivation for completeness. If the graph of color classes is a path, label the middle color $2$. Since the colors $1$, $2$ and the
colors $2$, $3$ are adjacent in this graph, the off-diagonal entries in \eqref{eq:path} are
positive, so $p,s\in\{1,2,3\}$ and $(q,r)\in\{(1,1),(1,2),(2,1)\}$. Interchanging the endpoint
colors $1$ and $3$ replaces $(p,s,q,r)$ by $(s,p,r,q)$. This allows $p\le s$ when $q=r=1$, and
allows $(q,r)=(1,2)$ otherwise. There are six possibilities of the first type and nine of the
second.

If it is a triangle, all off-diagonal entries are either $1$ or $2$, with at most one $2$ in each
row. By \eqref{eq:cycle}, the products of the entries on both sides are equal. This common product
cannot be $4$ or greater, since that would require at least four entries equal to $2$, whereas
there are at most three. Thus, the common product is either $1$ or $2$. If it is $1$, we obtain
$A_3$. If it is $2$, there are exactly two entries equal to $2$, one on each side of
\eqref{eq:cycle}. These entries are either opposite, giving $A_4$, or lie in the same column,
giving $C$, up to renaming colors.
\end{proof}

Three of the path matrices will be referred to by name:
\[
A_1=A(3,1,1,2)=\begin{pmatrix}0&3&0\\1&0&2\\0&1&2\end{pmatrix},\qquad
A_2=A(2,2,1,2)=\begin{pmatrix}1&2&0\\1&0&2\\0&2&1\end{pmatrix},
\]
\[
A_5=A(3,3,1,1)=\begin{pmatrix}0&3&0\\1&1&1\\0&3&0\end{pmatrix}.
\]

\section{Exclusion of the remaining candidates}\label{sec:exclusions}

We begin with a general observation about the path matrices \eqref{eq:path}, in which the colors
$1$ and $3$ are never adjacent.

\begin{lemma}\label{lem:outer2}
Let $T$ be a perfect $3$-coloring of $\GP(n,2)$ whose parameter matrix is a path matrix
\eqref{eq:path}. Then some vertex $a_i$ has color $2$.
\end{lemma}

\begin{proof}
Suppose that no vertex $a_i$ has color $2$. Then every $a_i$ has color $1$ or $3$, and since
$a_{13}=a_{31}=0$, the vertices $a_i$ and $a_{i+1}$ never have the colors $1$ and $3$. Hence all
$a_i$ have the same color $c\in\{1,3\}$. Every $a_i$ then has two
neighbors $a_{i-1},a_{i+1}$ of color $c$, so the color of its third neighbor $b_i$ is determined by
row $c$ of $A$ and is the same for every $i$. Thus at most two colors occur, contradicting the
assumption that $T$ is a perfect $3$-coloring.
\end{proof}

\begin{lemma}\label{lem:groupI}
None of the three matrices
\[
A(1,3,1,2)=\begin{pmatrix}2&1&0\\1&0&2\\0&3&0\end{pmatrix},\qquad
A(2,3,1,2)=\begin{pmatrix}1&2&0\\1&0&2\\0&3&0\end{pmatrix},\qquad
A(3,3,1,2)=\begin{pmatrix}0&3&0\\1&0&2\\0&3&0\end{pmatrix}
\]
occurs as a parameter matrix of a perfect $3$-coloring of $\GP(n,2)$.
\end{lemma}

\begin{proof}
First we show that no vertex $a_i$ has color $3$. Suppose $T(a_i)=3$. Then $a_{i-1}$, $a_{i+1}$,
and $b_i$ have color $2$, so
\[
\{T(a_{i+2}),T(b_{i+1})\}=\{T(a_{i-2}),T(b_{i-1})\}=\{T(b_{i-2}),T(b_{i+2})\}=\{1,3\}.
\]
If $T(a_{i+2})=3$, then its neighbor $b_{i+2}$ has color $2$, a contradiction. Hence
$T(a_{i+2})=1$ and $T(b_{i+1})=3$. But then the neighbor $b_{i-1}$ of $b_{i+1}$ has color $2$, a
contradiction.

Since color $3$ occurs, $T(b_i)=3$ for some $i$. Then $T(a_i)=2$, so $a_i$ must have two neighbors
of color $3$. Its neighbors $a_{i-1},a_{i+1}$ do not have color $3$ by the first part, and
$b_i$ is only one vertex. This is a contradiction.
\end{proof}

\begin{lemma}\label{lem:groupII}
Neither of the two matrices
\[
A(1,3,1,1)=\begin{pmatrix}2&1&0\\1&1&1\\0&3&0\end{pmatrix},\qquad
A(2,3,1,1)=\begin{pmatrix}1&2&0\\1&1&1\\0&3&0\end{pmatrix}
\]
occurs as a parameter matrix of a perfect $3$-coloring of $\GP(n,2)$.
\end{lemma}

\begin{proof}
\emph{Step 1: no vertex $a_i$ has color $3$.}
Suppose $T(a_i)=3$, so that $a_{i-1},a_{i+1},b_i$ have color $2$. Arguing as in
Lemma~\ref{lem:groupI}, we obtain
\[
\{T(a_{i+2}),T(b_{i+1})\}=\{1,2\},\qquad \{T(a_{i-2}),T(b_{i-1})\}=\{1,2\},\qquad
\{T(b_{i-2}),T(b_{i+2})\}=\{1,2\}.
\]
The adjacent vertices $b_{i-1}$ and $b_{i+1}$ cannot both have color $2$: otherwise $b_{i+1}$ would
have two neighbors of color $2$, namely $a_{i+1}$ and $b_{i-1}$. By symmetry, we may assume $T(b_{i+1})=1$, and therefore $T(a_{i+2})=2$. Then
$\{T(a_{i+3}),T(b_{i+2})\}=\{1,3\}$. Hence
$T(b_{i+2})=1$, $T(a_{i+3})=3$, and consequently $T(b_{i-2})=2$.
As $a_{i+3}$ has color $3$, its neighbor $b_{i+3}$ has color $2$.

Now $b_{i+1}$ has color $1$ and two neighbors of color $2$, namely $a_{i+1}$ and $b_{i+3}$. This is
impossible for $A(1,3,1,1)$, where a vertex of color $1$ has only one neighbor of color $2$. For
$A(2,3,1,1)$, the third neighbor $b_{i-1}$ of $b_{i+1}$ has color $1$, and hence $T(a_{i-2})=2$. But now $a_{i-2}$ has color $2$ and two neighbors of color $2$, namely $a_{i-1}$
and $b_{i-2}$, a contradiction.

\smallskip
\noindent\emph{Step 2.}
Since color $3$ occurs, Step~1 gives $T(b_j)=3$ for some $j$. Then $T(a_j)=2$. Since $a_j$ has exactly one neighbor of each
color and $b_j$ has color $3$, exactly one of $a_{j-1}$ and $a_{j+1}$ has color $2$. By symmetry, assume $T(a_{j+1})=2$. Hence
$T(a_{j+2})=1$ and $T(b_{j+1})=3$, which forces
$T(b_{j+2})=T(b_{j+3})=2$. Thus $a_{j+2}$ has color $1$ and two neighbors of color $2$, which is impossible for $A(1,3,1,1)$. For
$A(2,3,1,1)$, we get $T(a_{j+3})=1$,
forcing $T(a_{j+4})=2$ and hence $T(b_{j+4})=3$. But now $b_{j+2}$ has two neighbors of color $3$,
contradicting $a_{23}=1$.
\end{proof}

\begin{lemma}\label{lem:row102}
None of the four matrices
\[
A(1,1,1,2)=\begin{pmatrix}2&1&0\\1&0&2\\0&1&2\end{pmatrix},\qquad
A(1,2,1,2)=\begin{pmatrix}2&1&0\\1&0&2\\0&2&1\end{pmatrix},
\]
\[
A(2,1,1,2)=\begin{pmatrix}1&2&0\\1&0&2\\0&1&2\end{pmatrix},\qquad
A(3,2,1,2)=\begin{pmatrix}0&3&0\\1&0&2\\0&2&1\end{pmatrix}
\]
occurs as a parameter matrix of a perfect $3$-coloring of $\GP(n,2)$.
\end{lemma}

\begin{proof}
In all four matrices the second row is $(1,0,2)$ and $a_{13}=0$.
By Lemma~\ref{lem:outer2} there is a vertex $a_i$ of color $2$. Its neighbors
$a_{i-1},a_{i+1},b_i$ carry the colors $1,3,3$. By symmetry, there are two cases:
\begin{itemize}
\item[($\alpha$)] $T(b_i)=1$ and $T(a_{i-1})=T(a_{i+1})=3$,
\item[($\beta$)] $T(b_i)=3$, $T(a_{i+1})=1$, and $T(a_{i-1})=3$.
\end{itemize}
We treat the four matrices in turn.

\smallskip
\noindent\emph{The matrix $A(1,1,1,2)$.}
($\alpha$) Here $T(b_{i-2})=T(b_{i+2})=1$ and $T(a_{i+2})=T(b_{i+1})=3$. Then $b_{i+2}$, of color $1$, is adjacent to $a_{i+2}$, of color
$3$, contradicting $a_{13}=0$.
($\beta$) Here $T(a_{i+2})=T(b_{i+1})=1$ and $T(a_{i-2})=T(b_{i-1})=3$. Then $b_{i+1}$, of color $1$, is adjacent to
$b_{i-1}$, of color $3$, contradicting $a_{13}=0$.

\smallskip
\noindent\emph{The matrix $A(1,2,1,2)$.}
($\alpha$) Here $T(b_{i-2})=T(b_{i+2})=1$ and $\{T(a_{i+2}),T(b_{i+1})\}=\{2,3\}$. Since $b_{i+2}$ has color $1$, its neighbor $a_{i+2}$
does not have color $3$; so $T(a_{i+2})=2$ and $T(b_{i+1})=3$. By symmetry,
also $T(b_{i-1})=3$. Now $b_{i+1}$ has color $3$ and two neighbors of color $3$, namely
$a_{i+1}$ and $b_{i-1}$, contradicting $a_{33}=1$.
($\beta$) Here $T(a_{i+2})=T(b_{i+1})=1$, $\{T(a_{i-2}),T(b_{i-1})\}=\{2,3\}$, and
$\{T(b_{i-2}),T(b_{i+2})\}=\{2,3\}$. Since $b_{i+1}$ has color $1$, its neighbor $b_{i-1}$ has
color $2$, so $T(a_{i-2})=3$. Since $a_{i+2}$ has color $1$, its neighbor $b_{i+2}$ has color $2$,
so $T(b_{i-2})=3$. Now $b_{i-2}$ has color $3$ and two neighbors of color $3$, namely $b_i$ and
$a_{i-2}$, contradicting $a_{33}=1$.

\smallskip
\noindent\emph{The matrix $A(2,1,1,2)$.}
($\alpha$) Here $T(a_{i+2})=T(b_{i+1})=3$, $T(a_{i-2})=T(b_{i-1})=3$, and
$\{T(b_{i-2}),T(b_{i+2})\}=\{1,2\}$. Thus one of $b_{i-2}$ and $b_{i+2}$ has color $1$, but
its neighbor $a_{i-2}$ or $a_{i+2}$, respectively, has color $3$, contradicting $a_{13}=0$.
($\beta$) Here $\{T(a_{i+2}),T(b_{i+1})\}=\{1,2\}$, $T(a_{i-2})=T(b_{i-1})=3$, and
$T(b_{i-2})=T(b_{i+2})=3$. Since $b_{i+1}$ is adjacent to $b_{i-1}$, of color $3$, we get
$T(b_{i+1})=2$ and $T(a_{i+2})=1$. But $a_{i+2}$, of color $1$, is adjacent to $b_{i+2}$, of
color $3$, contradicting $a_{13}=0$.

\smallskip
\noindent\emph{The matrix $A(3,2,1,2)$.}
($\alpha$) Here $T(b_{i-2})=T(b_{i+2})=2$. Then $T(a_{i+2})=3$, and hence $T(b_{i+1})=2$. By symmetry,
also $T(b_{i-1})=2$. But $b_{i-1}$ and $b_{i+1}$ are adjacent, contradicting
$a_{22}=0$.
($\beta$) Here $T(a_{i+2})=T(b_{i+1})=2$. Then $T(b_{i-1})=T(b_{i+2})=3$, and hence
$T(a_{i-2})=T(b_{i-2})=2$. But $a_{i-2}$ and $b_{i-2}$ are adjacent, contradicting $a_{22}=0$.
\end{proof}

\begin{lemma}\label{lem:row111}
None of the three matrices
\[
A(1,1,1,1)=\begin{pmatrix}2&1&0\\1&1&1\\0&1&2\end{pmatrix},\qquad
A(1,2,1,1)=\begin{pmatrix}2&1&0\\1&1&1\\0&2&1\end{pmatrix},\qquad
A(2,2,1,1)=\begin{pmatrix}1&2&0\\1&1&1\\0&2&1\end{pmatrix}
\]
occurs as a parameter matrix of a perfect $3$-coloring of $\GP(n,2)$.
\end{lemma}

\begin{proof}
In all three matrices the second row is $(1,1,1)$ and $a_{13}=0$. Again by
Lemma~\ref{lem:outer2} there is a vertex $a_i$ of color $2$, and its neighbors
$a_{i-1},a_{i+1},b_i$ carry one color each. By the same argument as above, there are three cases:
\begin{itemize}
\item[($\alpha$)] $T(b_i)=2$, $T(a_{i+1})=1$, $T(a_{i-1})=3$,
\item[($\beta$)] $T(b_i)=1$, $T(a_{i+1})=3$, $T(a_{i-1})=2$,
\item[($\gamma$)] $T(b_i)=3$, $T(a_{i+1})=1$, $T(a_{i-1})=2$.
\end{itemize}

\smallskip
\noindent\emph{The matrix $A(1,1,1,1)$.}
($\alpha$) Here $T(a_{i+2})=T(b_{i+1})=1$ and $T(a_{i-2})=T(b_{i-1})=3$; but $b_{i+1}$ and $b_{i-1}$ are adjacent, contradicting $a_{13}=0$.
($\beta$) Here $T(b_{i-2})=T(b_{i+2})=1$ and $T(a_{i+2})=T(b_{i+1})=3$; but $b_{i+2}$ and $a_{i+2}$ are adjacent, contradicting $a_{13}=0$.
($\gamma$) Here $T(b_{i-2})=T(b_{i+2})=3$ and $T(a_{i+2})=T(b_{i+1})=1$; but $b_{i+2}$ and $a_{i+2}$ are adjacent, contradicting $a_{13}=0$.

\smallskip
\noindent\emph{The matrix $A(1,2,1,1)$.}
($\alpha$) Here $T(a_{i+2})=T(b_{i+1})=1$, $\{T(a_{i-2}),T(b_{i-1})\}=\{2,3\}$, and
$\{T(b_{i-2}),T(b_{i+2})\}=\{1,3\}$. Since $b_{i-1}$ is adjacent to $b_{i+1}$, of color $1$, we
get $T(b_{i-1})=2$ and $T(a_{i-2})=3$. Then $T(b_{i-2})=2$, contradicting $T(b_{i-2})\in\{1,3\}$.
($\beta$) Here $T(b_{i-2})=T(b_{i+2})=1$, $\{T(a_{i+2}),T(b_{i+1})\}=\{2,3\}$, and
$\{T(a_{i-2}),T(b_{i-1})\}=\{1,3\}$. Since $a_{i+2}$ is adjacent to $b_{i+2}$, of color $1$, we get
$T(a_{i+2})=2$ and $T(b_{i+1})=3$. Since $a_{i-2}$ is adjacent to $b_{i-2}$, of color $1$, we get
$T(a_{i-2})=1$ and $T(b_{i-1})=3$. Now $b_{i+1}$ has color $3$ and two neighbors of color $3$,
namely $a_{i+1}$ and $b_{i-1}$, contradicting $a_{33}=1$.
($\gamma$) Here $T(a_{i+2})=T(b_{i+1})=1$, $\{T(b_{i-2}),T(b_{i+2})\}=\{2,3\}$, and
$\{T(a_{i-2}),T(b_{i-1})\}=\{1,3\}$. Since $b_{i+2}$ is adjacent to $a_{i+2}$, of color $1$, we get
$T(b_{i+2})=2$ and $T(b_{i-2})=3$. Then $T(a_{i-2})=2$, contradicting $T(a_{i-2})\in\{1,3\}$.

\smallskip
\noindent\emph{The matrix $A(2,2,1,1)$.}
Interchanging the colors $1$ and $3$ leaves this matrix unchanged and maps case ($\gamma$) to case
($\beta$); so it suffices to treat ($\alpha$) and ($\beta$).
($\alpha$) Here $\{T(a_{i+2}),T(b_{i+1})\}=\{1,2\}$, $\{T(a_{i-2}),T(b_{i-1})\}=\{2,3\}$, and
$\{T(b_{i-2}),T(b_{i+2})\}=\{1,3\}$. If $T(b_{i+1})=2$, then $T(a_{i+2})=1$, and hence
$T(b_{i+2})=2$, a contradiction. If $T(b_{i+1})=1$, then $T(b_{i-1})=2$, hence $T(a_{i-2})=3$ and
$T(b_{i-2})=2$, a contradiction.
($\beta$) Here $\{T(b_{i-2}),T(b_{i+2})\}=\{1,2\}$, $\{T(a_{i+2}),T(b_{i+1})\}=\{2,3\}$, and
$\{T(a_{i-2}),T(b_{i-1})\}=\{1,3\}$. If $T(b_{i+1})=3$, then $T(b_{i-1})=2$, a contradiction. If
$T(b_{i+1})=2$, then $T(a_{i+2})=3$, hence $T(b_{i+2})=2$, $T(b_{i-2})=1$, and $T(a_{i-2})=2$, a
contradiction.
\end{proof}

\begin{lemma}\label{lem:C}
The triangle matrix $C$ does not occur as a parameter matrix of a perfect $3$-coloring of
$\GP(n,2)$.
\end{lemma}

\begin{proof}
First, some $a_i$ has color $3$. Otherwise all $a_i$ have colors $1$ and $2$, and since neither
color class contains an edge, the two colors alternate along the cycle $a_0a_1\cdots a_{n-1}$. This is
impossible if $n$ is odd. If $n$ is even, a vertex $a_i$ of color $1$ has two neighbors of color $2$,
contradicting $a_{12}=1$.

Let $T(a_i)=3$. Its neighbors $a_{i-1},a_{i+1},b_i$ carry one color each. Interchanging the colors $1$ and $2$
leaves $C$ unchanged. Using this interchange and symmetry, we may assume one of the following:
\begin{itemize}
\item[($\alpha$)] $T(b_i)=3$, $T(a_{i+1})=1$, $T(a_{i-1})=2$,
\item[($\beta$)] $T(b_i)=1$, $T(a_{i+1})=2$, $T(a_{i-1})=3$.
\end{itemize}

($\alpha$) Here $\{T(a_{i+2}),T(b_{i+1})\}=\{2,3\}$, $\{T(a_{i-2}),T(b_{i-1})\}=\{1,3\}$, and
$\{T(b_{i-2}),T(b_{i+2})\}=\{1,2\}$. If $T(a_{i+2})=2$, then $T(b_{i+2})=3$, a contradiction. If
$T(a_{i+2})=3$, then $T(b_{i+1})=2$, and hence $T(b_{i-1})=3$. Now $b_{i-1}$ has color $3$ and two neighbors of color $2$, namely $a_{i-1}$ and
$b_{i+1}$, contradicting $a_{32}=1$.

($\beta$) Here $\{T(b_{i-2}),T(b_{i+2})\}=\{2,3\}$, $\{T(a_{i+2}),T(b_{i+1})\}=\{1,3\}$, and
$\{T(a_{i-2}),T(b_{i-1})\}=\{1,2\}$. If $T(a_{i+2})=1$, then $T(b_{i+2})=3$; now $b_{i+2}$ has color $3$ and two
neighbors of color $1$, namely $b_i$ and $a_{i+2}$, contradicting $a_{31}=1$. If $T(a_{i+2})=3$,
then $T(b_{i+1})=1$, and hence $T(b_{i-1})=3$, a contradiction.
\end{proof}

\begin{proposition}\label{prop:list}
Every perfect $3$-coloring of $\GP(n,2)$ has parameter matrix $A_1$, $A_2$, $A_3$, $A_4$, or
$A_5$, up to renaming colors.
\end{proposition}

\begin{proof}
This is a direct consequence of Lemmas~\ref{lem:candidates} and \ref{lem:groupI}--\ref{lem:C}.
\end{proof}

\begin{remark}
In particular, Lemmas~\ref{lem:groupII}, \ref{lem:row102} and \ref{lem:row111} contain detailed
proofs, valid for every $n\ge 5$, of the five exclusions in \cite[Lemma~3.4]{AKS}, which were used
there for $n=6$ and $n=7$.
\end{remark}

\section{The remaining matrices}\label{sec:explicit}

\begin{proposition}\label{prop:words}
The following mappings $T\colon V(\GP(n,2))\to\{1,2,3\}$ are perfect $3$-colorings with the
indicated parameter matrices. In each case $i$ runs over $0,1,\dots,m-1$.
\begin{itemize}
\item[$A_1$:] Let $n=5m$. Put $T(a_{5i})=1$ and
$T(a_{5i+1})=T(a_{5i+4})=T(b_{5i})=2$, and give all other vertices color $3$.
\item[$A_2$:] Let $n=5m$. Put $T(a_{5i})=T(b_{5i})=1$ and
$T(a_{5i+1})=T(a_{5i+4})=T(b_{5i+2})=T(b_{5i+3})=2$, and give all other vertices color $3$.
\item[$A_3$:] Let $n=3m$. Put $T(a_{3i})=T(b_{3i})=1$, $T(a_{3i+1})=T(b_{3i+1})=2$, and
$T(a_{3i+2})=T(b_{3i+2})=3$.
\item[$A_4$:] Let $n=3m$. Put $T(a_{3i})=T(b_{3i+2})=1$, $T(a_{3i+1})=T(b_{3i+1})=2$, and
$T(a_{3i+2})=T(b_{3i})=3$.
\item[$A_5$:] Let $n=10m$. Put
$T(a_{10i})=T(a_{10i+3})=T(b_{10i+6})=T(b_{10i+7})=1$ and
$T(a_{10i+5})=T(a_{10i+8})=T(b_{10i+1})=T(b_{10i+2})=3$, and give all other vertices color $2$.
\end{itemize}
\end{proposition}

\begin{proof}
It is trivial.
\end{proof}

\begin{proposition}\label{prop:div}
The matrices $A_1,A_2$ occur exactly when $5\mid n$, and $A_3,A_4$ occur exactly when $3\mid n$.
\end{proposition}

\begin{proof}
Equation~\eqref{eq:balance}, together with $v_1+v_2+v_3=2n$, gives
\[
\begin{array}{c|ccc}
\text{Matrix} & v_1 & v_2 & v_3\\ \hline
A_1 & \dfrac{n}{5} & \dfrac{3n}{5} & \dfrac{6n}{5}\\[8pt]
A_2 & \dfrac{2n}{5} & \dfrac{4n}{5} & \dfrac{4n}{5}\\[8pt]
A_3,A_4 & \dfrac{2n}{3} & \dfrac{2n}{3} & \dfrac{2n}{3}
\end{array}
\]
Since $v_1,v_2,v_3$ are integers, we obtain $5\mid n$ for $A_1,A_2$ and $3\mid n$ for $A_3,A_4$.
The converse follows from Proposition~\ref{prop:words}.
\end{proof}

\begin{remark}
For $A_3$ and $A_4$, every vertex has exactly one neighbor of color $2$. Hence the color class $2$
is a total perfect code of $\GP(n,2)$. By \cite[Theorem~1.2]{WZ}, all total perfect codes of
$\GP(n,2)$ are precisely $C_j=\{a_{3i+j},b_{3i+j}: 0\le i<n/3\}$, $j=0,1,2$, whenever $3\mid n$, and there are none otherwise.
\end{remark}

For $A_5$, Equation~\eqref{eq:balance}, together with $v_1+v_2+v_3=2n$, gives
$v_1=v_3=\frac{2n}{5}$ and $v_2=\frac{6n}{5}$, so $5\mid n$ is necessary.

\begin{proposition}\label{prop:A5}
There is a perfect $3$-coloring of $\GP(n,2)$ with parameter matrix $A_5$ if and only if
$10\mid n$.
\end{proposition}

\begin{proof}
Not every $a_i$ can have color $2$: otherwise each $a_i$
would have the two color-$2$ neighbors $a_{i-1}$ and $a_{i+1}$. By symmetry, and since interchanging colors $1$ and $3$ leaves $A_5$
unchanged, we may assume $T(a_0)=1$. Then
$T(a_{-1})=T(a_1)=T(b_0)=2$. Since $b_0$ has one neighbor of each color, the colors on $b_{-2},b_2$ are $2$ and $3$. By symmetry, we may take
$T(b_2)=3$ and $T(b_{-2})=2$. Since $b_2$ has color $3$, the vertices $a_2,b_4$ have color $2$. Now $T(b_1)=3$ and
$T(a_3)=1$. Thus $T(a_4)=T(b_3)=2$. At $b_3$, the
other known neighbors $a_3,b_1$ have colors $1$ and $3$, so $T(b_5)=2$. At $a_4$, the neighbors
$a_3,b_4$ have colors $1$ and $2$, so $T(a_5)=3$. Finally,
$T(b_6)=T(b_7)=1$.

At index $5$ the initial configuration is repeated, with colors $1$ and $3$ interchanged: $a_5$ has
color $3$, $b_5$ has color $2$, $b_7$ has color $1$, and $b_3$ has color $2$. Repeating the
argument determines every subsequent color. Each shift by five interchanges colors $1$ and $3$ and
fixes color $2$; every shift by ten restores all colors. The forced coloring is exactly the one for
$A_5$ in Proposition~\ref{prop:words}.

In particular, for $j\ge 0$ the forced color of $a_j$ is $1$ exactly when $j=10i$ or $j=10i+3$
for some $i\ge 0$. Since $a_n=a_0$ and $a_{n+3}=a_3$ have color $1$, both $n$ and $n+3$ are of this
form, which is possible only if $10\mid n$. Conversely, Proposition~\ref{prop:words} proves
existence.
\end{proof}

\needspace{8\baselineskip}
\section{The classification}\label{sec:main}

We now put the previous results together.

\begin{theorem}\label{thm:main}
Let $n\ge 5$. Up to renaming the colors, every parameter matrix of a perfect $3$-coloring of
$\GP(n,2)$ is one of $A_1,\dots,A_5$, and these matrices occur exactly under the following
conditions:
\[
\begin{array}{c|c}
\text{Matrix} & \text{Necessary and sufficient condition}\\ \hline
A_1,A_2 & 5\mid n\\
A_3,A_4 & 3\mid n\\
A_5 & 10\mid n.
\end{array}
\]
\end{theorem}

\begin{proof}
This is a direct consequence of Propositions~\ref{prop:list}, \ref{prop:div} and \ref{prop:A5}.
\end{proof}

\begin{remark}
Theorem~\ref{thm:main} agrees with the cases treated earlier: $n=5,6,7$ in \cite{AKS}, $n=8$ in
\cite{AS}, and $n=5$ and $n=10$ in \cite{AQ}.
\end{remark}

Let $C$ be a nonempty set of vertices of a connected graph. If coloring each vertex by its distance
from $C$ gives a perfect coloring, then $C$ is a \emph{completely regular code}, and the largest
distance of a vertex from $C$ is its \emph{covering radius} \cite{KP25}.

The perfect $3$-colorings whose graph of color classes is a path are exactly the distance
partitions of completely regular codes of covering radius two. Thus the path matrices $A(p,s,q,r)$
in \eqref{eq:path} are the quotient matrices of such codes; in the usual notation, the code formed
by color class $1$ has intersection array $\{p,r;q,s\}$ \cite{BRZ19}.

\begin{corollary}\label{cor:crc}
Let $n\ge 5$. The graph $\GP(n,2)$ contains a completely regular code of covering radius two if and
only if $5\mid n$. For such a code $C$, the parameter matrix of the coloring by distance from $C$ is
$A_1$, $A_2$, or $A_5$, possibly with the order of the colors reversed, and $A_5$ occurs if and only
if $10\mid n$.
\end{corollary}

\end{document}